\documentclass[11pt,leqno]{amsart}
\usepackage{amssymb,verbatim,enumerate,ifthen}
\usepackage[mathscr]{eucal}
\usepackage[utf8]{inputenc}
\usepackage[T1]{fontenc}
\def\N{\mathbb{N}}
\def\R{\mathbb{R}}

\newtheorem{theorem}{Theorem}[section]
\newtheorem*{theorem*}{Theorem}
\def\Thm#1#2{\ifthenelse{\equal{#1}{*}}{\begin{theorem*}#2\end{theorem*}}
             {\begin{theorem}\label{T#1}#2\end{theorem}}}
\newtheorem{Atheorem}{Theorem}

\newtheorem{proposition}[theorem]{Proposition}
\newtheorem*{proposition*}{Proposition}
\def\Prp#1#2{\ifthenelse{\equal{#1}{*}}{\begin{proposition*}#2\end{proposition*}}
{\begin{proposition}\label{P#1}#2\end{proposition}}}
\def\prp#1{Proposition~\ref{P#1}}

\newtheorem{corollary}[theorem]{Corollary}
\newtheorem*{corollary*}{Corollary}
\def\Cor#1#2{\ifthenelse{\equal{#1}{*}}{\begin{corollary*}#2\end{corollary*}}
             {\begin{corollary}\label{C#1}#2\end{corollary}}}

\newtheorem{lemma}[theorem]{Lemma}
\newtheorem*{lemma*}{Lemma}
\def\Lem#1#2{\ifthenelse{\equal{#1}{*}}{\begin{lemma*}#2\end{lemma*}}
             {\begin{lemma}\label{L#1}#2\end{lemma}}}

\theoremstyle{definition}
\newtheorem{remark}[theorem]{Remark}
\newtheorem*{remark*}{Remark}
\def\Rem#1#2{\ifthenelse{\equal{#1}{*}}{\begin{remark}\rm #2\end{remark}}
             {\begin{remark}\label{R#1}\rm #2\end{remark}}}

\newtheorem{example}[theorem]{Example}
\newtheorem*{example*}{Example}
\def\Exa#1#2{\ifthenelse{\equal{#1}{*}}{\begin{example*}\rm #2\end{example*}}
             {\begin{example}\label{Ex#1}\rm #2\end{example}}}

\def\eq#1{{\rm(\ref{E#1})}}
\def\Eq#1#2{\ifthenelse{\equal{#1}{*}}
  {\begin{equation*}\begin{aligned}#2\end{aligned}\end{equation*}}
  {\begin{equation}\begin{aligned}\label{E#1}#2\end{aligned}\end{equation}}}

\begin{document}
\begin{flushright}
\end{flushright}
\vspace{5mm}

\date{\today}

\title{On Characterizations of Convex and \\
Approximately Subadditive Sequences}

\author[A. R. Goswami]{Angshuman R. Goswami}
\address[A. R. Goswami]{Department of Mathematics, University of Pannonia, 
H-8200 Veszprem, Hungary}
\email{goswami.angshuman.robin@mik.uni-pannon.hu}
\thanks{The author’s research was supported by the EKÖP Scholarship (2024-2.1.1-EKÖP-2024-00025/29), funded by the Research Fellowship Programme of the Ministry of Culture and Innovation, Government of Hungary.}

\begin{abstract}
A sequence $\Big(u_n\Big)_{n=0}^{\infty}$ is said to be convex if it satisfies the following inequality
\Eq{*}{
2u_n\leq u_{n-1}+u_{n+1}\qquad \mbox{for all}\qquad n\in\N.
}
We present several characterizations of convex sequences and demonstrate that such sequences can be locally interpolated by quadratic polynomials. Furthermore, the converse assertion of this statement is also established.\\

On the other hand,  a sequence $\Big(u_n\Big)_{n=1}^{\infty}$ is called  approximately subadditive if for a fixed $\epsilon>0$ and any partition $n_1,\cdots,n_k$ of $n\in\N$;  the following discrete functional inequality holds true
\Eq{*}{
u_n\leq u_{n_1}+\cdots+ u_{n_k}+\varepsilon.
}
We show Ulam's type stability result for such sequences. 
We prove that an approximately subadditive sequence can be expressed as the algebraic summation of an ordinary subadditive and a non-negative sequence bounded above by $\varepsilon.$ \\

A proposition portraying the linkage between the convex and subadditive sequences under minimal assumption is also included.\\

The motivation, research background, important notions, and terminologies are discussed in the introduction section.
\end{abstract}

\maketitle

\section*{Introduction} 
Throughout this paper $\N$, $\R$, and $\R_+$ denote the sets of natural, real, and positive numbers respectively. This paper primarily aims to introduce multiple characterizations of convex sequences and to present a decomposition result concerning approximately subadditive sequences. \\

A sequence $\Big(u_n\Big)_{n=0}^{\infty}$ is said to be \textit{sequentially convex} (or a \textit{convex sequence}) if for all $n\in\N$, it satisfies the following functional inequality
\Eq{2000}{
2u_n\leq u_{n-1}+u_{n+1}.
}
If the converse of the above inequality holds, $\Big(u_n\Big)_{n=0}^{\infty}$ would be termed as a \textit{concave sequence}. Arithmetic, geometric, Fibonacci, and many other notable sequences can be analyzed within the framework of convex sequences.
\\

Although the first mention of the convex sequence is not very clear; based on limited evidence, most mathematicians believe that the terminology of sequential convexity first appeared in the book of \cite{Mitrinovicc}. Some of the early works in this direction can be found in the papers \cite{Essen, pecaric, GauSte}.\\

Since then, mathematicians have also explored various aspects of sequential convexity, including discrete analogous of the Hermite-Hadamard inequality, linkage with difference equations, applications in numerical estimation and trigonometric functions etc. Relevant findings can be found in the works of \cite{Latreuch, Debnath, GAB_PAL}.\\

In recent times, researchers have studied several new versions of convex sequences.  For example, investigations into higher-ordered convex,  relatively convex, symmetrized convex, approximately convex, and $\alpha$-convex sequences have been carried out. The findings enhance the understanding of functional inequalities in discrete settings. The studies in \cite{Milovanovic, Jim\'enez, Mercer, HDCS, Farissi, Abdellaoui, Goswamii, Goswami, Krasniqi} offer insight into these topics.
\\

In the first section of this paper, we provide some characterizations of convex sequences. We show that for a convex sequence $\Big(u_n\Big)_{n=0}^{\infty}$, there exists an underlying monotone sequence $\Big(v_n\Big)_{n=0}^{\infty}$ that tracks the discrete slope between successive terms. Also, we prove that any convex sequence can be locally interpolated by a quadratic polynomial. More details regarding such generalizations, characterizations, decompositions, and other such related research of different function and sequence classes can be found in the books \cite{Kuczma, Hardy, Schechter}.\\

For any chosen $n\in\N$, a finite collection $\{n_1,\cdots n_k\}\subset\N$ is called a \textit{partition} of $n$ if 
\newline $n=n_1+\cdots+n_k$ holds. 
Partitions play a significant role in the field of number theory. Using the concept of partition, we define sequential subadditivity.\\

Let $\varepsilon>0$ and $n\in\N$ be arbitrary. A sequence $\Big(u_n\Big)_{n=1}^{\infty}$ is said to be \textit{approximately subadditive} if for any partition $n_1,\cdots,n_k $ of $n$, the following discrete functional inequality holds true
\Eq{0}{
u_n\leq u_{n_1}+\cdots+u_{n_k}+\varepsilon.
}
On the other hand, if the sequence $\Big(u_n\Big)_{n=1}^{\infty}$ satisfies the above inequality without the $\epsilon$; in such case we said it as an ordinary \textit{subadditive sequence}.\\

Discrete subadditivity is at the center of many important mathematical results. One can look into Fekete's subadditive lemma \cite{Fekete}, Kingman's subadditive Ergodic theorem \cite{Kingman}, Hammersley's subadditive theorem \cite{Hammersley} etc. In recent years, researchers investigated several new variations of sequential subadditivities. Alternative proofs for some of the well-known results concerning subadditive sequences are also provided. To explore this further, see the results reported in \cite{Zoltan, Angshuman}. Based on all these research several applications in optimal transport, machine learning, information theory, and mathematical modelling are also proposed.\\

\textit{'Ulam-type stability'} studies when an approximate solution to a functional equation implies the existence of a true solution nearby. Originating from Ulam’s 1940 problem and Hyers' subsequent answer, it formalizes that small deviations in functional relations do not drastically alter the solution space. This notion has been extended to various mathematical settings, including differential equations, group homomorphisms, and convexity. Stability results of this type are fundamental in understanding the robustness and rigidity of mathematical structures. Additional details are available in the classical papers \cite{Hyers, Ulam}. \\

However, stability analysis is relatively a new concept in sequence settings. In the second section of this paper, we propose a Ulam-type stability theorem for subadditive sequences.\\

It is evident that if $\Big(v_n\Big)_{n=1}^{\infty}$ possesses sequential subadditivity and $\Big(w_n\Big)_{n=1}^{\infty}$ is a non-negative sequence bounded above by $\varepsilon$, then the derived sequence $\Big(u_n\Big)_{n=1}^{\infty}: =\Big(v_n\Big)_{n=1}^{\infty}+\Big(w_n\Big)_{n=1}^{\infty}$ is a approximately subadditive majorant of $\Big(v_n\Big)_{n=1}^{\infty}$. Also, we estalished that if a sequence $\Big(u_n\Big)_{n=1}^{\infty}$ satisfies the inequality \eq{0}, then it can be expressed as the algebraic summation of two sequences $\Big(v_n\Big)_{n=1}^{\infty}$ and $\Big(w_n\Big)_{n=1}^{\infty}$. Where $\Big(v_n\Big)_{n=1}^{\infty}$ is a subadditive minorant of $\Big(u_n\Big)_{n=1}^{\infty}$; while $\Big(w_n\Big)_{n=1}^{\infty}$ is a non-negative sequence bounded above by $\varepsilon$.\\

It is a interesting observation that the ordering of convex sequences starts from index $0$, whereas for approximately subadditive (or ordinary subadditive) sequences, we begin with index $1$. This distinction arises because, in subadditive sequence classes, the ordering is crucial. For instance, in order for \eq{0} to hold, it is important first to consider all partitions of the positive integer $n$, followed by the corresponding sequential values at the partitioning points. In contrast, each term in a convex sequence depends only on the average of its two neighbouring terms, making the global ordering  redundant. These distinctions are discussed in detail in \cite{Rob}.\\

Now, we start our investigation from convex sequences.

\section{Characterization of Convex Sequences}
We begin with a fundamental fractional inequality that will play a central role in several subsequent results. This inequality is also mentioned in one of our recently submitted papers. However, for readability purposes, we state the statement and propose a shorter proof.
\Lem{22}{Let $a_1,\cdots,a_n\in\R$ and $b_1,\cdots,b_n\in\R_+$, then the following discrete functional inequality is satisfied
\Eq{1}{
\min\Bigg(\dfrac{a_1}{b_1},\cdots,\dfrac{a_n}{b_n}
\Bigg)\leq 
\dfrac{a_1+\cdots+a_n}{b_1+\cdots+b_n}
\leq\max\Bigg(\dfrac{a_1}{b_1},\cdots,\dfrac{a_n}{b_n}
\Bigg).
}
}
\begin{proof}
The expression $\dfrac{a_1+\cdots+a_n}{b_1+\cdots+b_n}$ can be re-written as the following convex combination
\Eq{*}{
\dfrac{b_1}{b_1+\cdots+b_n}\Bigg(\dfrac{a_1}{b_1}\Bigg)+\cdots+\dfrac{b_n}{b_1+\cdots+b_n}\Bigg(\dfrac{a_n}{b_n}\Bigg).
}
Hence by mean property the inequality \eq{1} is obvious.
\end{proof}
\Thm{11}{Let $\Big(u_n\Big)_{n=0}^{\infty}$ be a real valued sequence. Then the following conditions are equivalent to each other
\begin{enumerate}[(i)]
\item $\Big(u_n\Big)_{n=0}^{\infty}$ is convex.\\
\item For all $n_1,n_2,n_3\in \N\cup\{0\}$ with $n_1<n_2<n_3$, it satisfies
 \Eq{2}{
 \dfrac{u_{n_2}-u_{n_1}}{n_2-n_1}\leq \dfrac{u_{n_3}-u_{n_2}}{n_3-n_2}.
 }
\item There exists a monotone (increasing) sequence $\Big(v_n\Big)_{n=0}^{\infty}$ such that for all $m,n\in \N\cup\{0\}$,
\Eq{4}{
u_n-u_m\leq v_n(n-m) \qquad (m,n\in \N\cup\{0\}).
} 
\end{enumerate}
}
\begin{proof} {\it(i)$\rightarrow$(ii):} Assume that $\Big(u_n\Big)_{n=0}^{\infty}$ possesses sequential convexity and let $n_1,n_2$ and $n_3\in\N\cup\{0\}$ with $n_1<n_2<n_3$. Using the inequality \eq{1}, we proceed as follows
\Eq{100}{
\dfrac{u_{n_2}-u_{n_1}}{n_2-n_1}&=\dfrac{u_{n_2}-u_{n_2-1}+\cdots +u_{n_1+1}-u_{n_1}}{1+\cdots+1}\\
&\leq \max{\Big(u_{n_2}-u_{n_2-1},\cdots, u_{n_1+1}-u_{n_1}\Big)}\\
&=u_{n_2}-u_{n_2-1}.
}
Similarly,
\Eq{200}{
\dfrac{u_{n_3}-u_{n_2}}{n_3-n_2}&=\dfrac{u_{n_3}-u_{n_3-1}+\cdots +u_{n_2+1}-u_{n_2}}{1+\cdots+1}\\
&\geq \min{\Big(u_{n_3}-u_{n_3-1},\cdots, u_{n_2+1}-u_{n_2}\Big)}\\
&=u_{n_2+1}-u_{n_2}.
}
Convexity of the sequence $\Big(u_n\Big)_{n=0}^{\infty}$ implies $u_{n_2}-u_{n_2-1}\leq u_{n_2+1}-u_{n_2}$. This along with \eq{100} and \eq{200} establishes \eq{2} and completes the assertion. 
\\

{\it(ii)$\rightarrow$(iii):}
Assume that {\it(ii)} holds. We define the sequence 
$\Big(v_n\Big)_{n=0}^{\infty}$ as follows
\Eq{*}{
 v_n:=\underset{n\leq n_1<n_2}{\mathrm{\inf}}\Bigg(\dfrac{u_{n_2}-u_{n_1}}{n_2-n_1}\Bigg) \qquad(n_1, n_2\in \N\cup\{0\}).
}
In view of condition {\it(ii)}, for all $n_1<n_2<n_3$ in $\N\cup\{0\}$, we can conclude
\Eq{6}{
\dfrac{u_{n_2}-u_{n_1}}{n_2-n_1}\leq v_{n_2}\leq \dfrac{u_{n_3}-u_{n_2}}{n_3-n_2}.
 }
From the left-most inequality of \eq{6}, we get
\Eq{57}{
u_{n_2}-u_{n_1}\leq v_{n_2} (n_2-n_1) \qquad (n_1,n_2\in \N\cup\{0\}\quad\mbox{with}\quad n_1<n_2).
}
Similarly, from the right-most inequality of \eq{6} (replacing 
$n_3$ with $n_1$), it follows that
\Eq{*}{
u_{n_2}-u_{n_1}\leq v_{n_2} (n_2-n_1) \qquad (n_1,n_2 \in \N \quad\mbox{with}\quad n_2<n_1).
}
Also, by definition of $\Big(v_n\Big)_{n=0}^{\infty}$,  it is obvious that
$v_0\leq u_1-u_0.$ This together with the above inequality yields 
\Eq{58}{
u_{n_2}-u_{n_1}\leq v_{n_2} (n_2-n_1)  \qquad (n_1,n_2 \in \N\cup\{0\} \quad\mbox{with}\quad n_2<n_1).
}

The combined \eq{57} and \eq{58} can be represented as the following generalized form of inequality 
\Eq{*}{
u_n-u_m\leq v_n(n-m) \qquad \quad (\,\,\mbox{for all}\quad m,n\in \N\cup\{0\}\,\,).
}
This is the inequality \eq{4} which was needed to be established.
\\

Now to show the monotonicity of the sequence $\Big(v_n\Big)_{n=0}^{\infty}$, we assume $m,n\in\N\cup\{0\}$ with $m<n$. Then the assertion $(iii)$ provides the following two discrete inequalities
\Eq{*}{
u_n-u_m\leq v_n(n-m) \qquad \mbox{and}\qquad u_m-u_n\leq v_m(m-n).
}
Summing up these two inequalities side by side, we arrive at
\Eq{*}{
0 \leq (v_n-v_m)(n-m) \qquad m,n\in\N\cup\{0\} \quad \mbox{with} \quad m<n.
}
This implies $v_m\leq v_n$. Since $m,n\in\N\cup\{0\}$ are arbitrarily chosen, hence $\Big(v_n\Big)_{n=0}^{\infty}$ possesses monotonicity.\\

{\it(iii)$\rightarrow$(i):} Now we assume that the condition (iii) holds. Let $n\in\N$ be arbitrary.
By replacing $m$ with $n-1$ and $n+1$ respectively in the inequality \eq{4}, we obtain the following two inequalities
\Eq{*}{
u_n-u_{n-1}\leq v_n(n-(n-1))\qquad\mbox{and}\qquad u_n-u_{n+1}\leq v_n(n-(n+1)).
}
Summing up these two inequalities side by side, we arrive at
\eq{2000}.
This implies convexity of the sequence $\Big(u_n\Big)_{n=0}^{\infty}$ and completes the proof.
\end{proof}
The Lagrange polynomial of degree $n$ associated with the sequence $\Big(u_n\Big)_{n=0}^{\infty}$ can be expressed as follows
\Eq{*}{
P_n(x):=\sum_{i=0}^n\Bigg(\prod_{i\neq j}\dfrac{x-x_j}{i-j}u_i\Bigg)\qquad (n\in\N).
}
For a convex sequence $\Big(u_n\Big)_{n=0}^{\infty}$, the interpolating  Lagrange polynomial does not possess convexity.  For instance (0,-1,1,3) is a convex sequence. But the  corresponding Lagrange polynomial
\Eq{*}{
P_3(x)=-\dfrac{1}{2}x^3+3x^2-\dfrac{7}{2}x, \qquad x\in[0,3]}
is neither convex nor concave. However, in the next proposition, we will see that any convex sequence can be locally interpolated by a spline of degree 2. \\ 

To construct the result, corresponding to the sequence 
$\Big(u_n\Big)_{n=0}^{\infty}$, for each $n\in\N$, we define Lagrange polynomials of degree $2$ in $[n-1,n+1]$ as follows

\Eq{20}{
P_{_2}^{^{n}}(x):=\bigg(\dfrac{u_{n-1}+u_{n+1}}{2}-u_n\bigg)x^2+\Bigg(\dfrac{u_{n+1}-u_{n-1}}{2}-2\bigg(\dfrac{u_{n-1}+u_{n+1}}{2}-u_n\bigg) n\Bigg)x\\
\qquad \qquad \qquad\,\,\,
+\Bigg(\bigg(\dfrac{u_{n-1}+u_{n+1}}{2}-u_n\bigg)n^2-\bigg(\dfrac{u_{n+1}-u_{n-1}}{2}\bigg)n+u_n \Bigg).
}

\Prp{1}{A sequence $\Big(u_n\Big)_{n=0}^{\infty}$ is convex if and only if the quadratic polynomial $P_{_2}^{^{n}}(x)$ defined in \eq{20} is convex for each $n\in\N$.}
\begin{proof}
To prove the proposition, first we assume that $\Big(u_n\Big)_{n=0}^{\infty}$ is convex. Using convexity of the sequence, it can be easily observed that
\Eq{3}{\Big(P_{_2}^{^{n}}\Big)^{''}=u_{n-1}-2u_n+u_{n+1}\geq 0 \qquad \mbox{for all} \quad n\in\N\,\,;}
where $\Big(P_{_2}^{^{n}}\Big)^{''}$ denotes the second derivative of $P_{_2}^{^{n}}$.
Hence, $P_{_2}^{^{n}}$ is convex on the interval $[n-1,n+1]$ for each $n\in\N.$\\

To prove the converse part, we assume that for each $n\in\N$ the polynomial $P_{_2}^{^{n}}$ is convex. In other words, \eq{3} is satisfied. The right-most inequality of \eq{3} is just re-arranged form of \eq{2000}.  This establishes that the sequence $\Big(u_n\Big)_{n=0}^{\infty}$ is convex and completes the proof.
\end{proof}

\section{Decomposition of Approximately Subadditive Sequences}
In this section, we discuss the decomposition of approximately subadditive sequences. The proof of the proposed result extensively utilizes proposition 1.2 of our paper \cite{Rob}.\\

\Prp{77}{The sequence $\Big(u_n\Big)_{n=1}^{\infty}$ is approximately subadditive if and only if it can be expressed as the algebraic summation of a subadditive minorant 
$\Big(v_n\Big)_{n=1}^{\infty}$ and a non-negative sequence $\Big(w_n\Big)_{n=1}^{\infty}$ which is bounded above by $\varepsilon$.\\ 
}
\begin{proof}
First, we consider $n\in\N$ to be arbitrary and $n_1,\cdots,n_k$ to be its arbitrary partition. Then using the subadditivity and the bounds of the respective sequence $\Big(v_n\Big)_{n=1}^{\infty}$ and $\Big(w_n\Big)_{n=1}^{\infty}$, we can compute the following
\Eq{123}{
v_n+w_n\leq v_{n_1}+\cdots+v_{n_k}+\varepsilon.
}
We define $\Big(u_n\Big)_{n=1}^{\infty}=\Big(v_n\Big)_{n=1}^{\infty}+\Big(w_n\Big)_{n=1}^{\infty}.$ Here, the non-negativity of the sequence $\Big(w_n\Big)_{n=1}^{\infty}$ implies $v_n\leq u_n$ for all $n\in\N.$ Hence, the inequality \eq{123} can be extended as follows
\Eq{*}{
u_n\leq v_{n_1}+\cdots+v_{n_k}+\varepsilon\leq 
u_{n_1}+\cdots+u_{n_k}+\varepsilon.
}
This shows that the sequence  $\Big(u_n\Big)_{n=1}^{\infty}$ is approximately subadditive.\\

To prove the converse part, we assume that $\Big(u_n\Big)_{n=1}^{\infty}$ be approximately subadditive. In other words, for the sequence $\Big(u_n\Big)_{n=1}^{\infty}$, the inequality \eq{0} is satisfied. Now, we assume $n\in\N$ and $n_1,\cdots,n_k$ be any arbitrary partition of it. We construct the sequence $\Big(v_n\Big)_{n=1}^{\infty}$ as follows
\Eq{701}{
v_n:=\min\Big\{u_{n_1}+\cdots+u_{n_k}\,\,\, \Big|\,\,\, n_1,\cdots,n_k\in\N  \quad \mbox{satisfying}\quad n_1+\cdots+n_k=n\Big\}
.}
Clearly $v_n\leq u_n$ holds for all $n\in\N$. We only need to show that the sequence $\Big(v_n\Big)_{n=1}^{\infty}$ is subadditive. \\

We consider $m,n\in\N$ and have their respective partitions such that $m=m_1+\cdots +m_l$ and $n=n_1+\cdots+n_k$ satisfying the following two discrete functional equalities
\Eq{33}{
v_m=u_{m_1}+\cdots+u_{m_l}\quad\mbox{and}\quad
v_n=u_{n_1}+\cdots+u_{n_k}.}
The combined partitions of $m$ and $n$ provide a partition for $m+n$ as well. This can be represented as follows   
$$m+n= m_1+\cdots m_l+ n_1+\cdots+ n_k.$$

From the construction of the sequence  $\Big(v_n\Big)_{n=1}^{\infty}$ (inequality \eq{701}) and using \eq{33}, we can conclude the following inequality
\Eq{*}{
v_{m+n}\leq u_{m_1}+\cdots+u_{m_l}+ u_{n_1}+\cdots+u_{n_k}= v_m+v_n.
}
This yields that sequence $\Big(v_n\Big)_{n=1}^{\infty}$ is subadditive.\\

Now we define $\Big(w_n\Big)_{n=1}^{\infty}:=\Big(u_n\Big)_{n=1}^{\infty}-\Big(v_n\Big)_{n=1}^{\infty}$. This ensures the non-negativity of the sequence 
$\Big(w_n\Big)_{n=1}^{\infty}.$ To determine its upper bound, first, we consider any point $w_n$ from the sequence $\Big(w_n\Big)_{n=1}^{\infty}$. For the $n\in\N$, there must exists a partition 
$\{n_1,\cdots,n_k\}$ that satisfies the second equality of \eq{33}. Using this together with \eq{701}, we can  proceed as follows
\Eq{*}{
w_n=u_n-v_n\leq (u_{n_1}+\cdots+u_{n_k}+\varepsilon)-(u_{n_1}+\cdots+u_{n_k})=\varepsilon.
}
This shows that the sequence $\Big(w_n\Big)_{n=1}^{\infty}$ is bounded and completes the proof.
\end{proof}
The following proposition is self-verifiable. Hence the proof is left for the reader.
\Prp{78}{If all the elements of the decreasing sequence $\Big(v_n\Big)_{n=1}^{\infty}$ are non-negative, then it also possesses sequential subadditivity.}
The next proposition establishes a connection between the concave and subadditive sequences
\Prp{79}{
If all the elements of concave and monotone(increasing) sequence $\Big(u_n\Big)_{n=0}^{\infty}$ are non-negative, then the corresponding sequence 
$\Big(u_n-u_{n-1}\Big)_{n=1}^{\infty}$ is subadditive.
Conversely,  if the sequence 
\newline
$\Big(u_n-u_{n-1}\Big)_{n=1}^{\infty}$ is decreasing, then the corresponding sequence $\Big(u_n\Big)_{n=0}^{\infty}$ is concave.
}
\begin{proof}
The concavity of the sequence $\Big(u_n\Big)_{n=0}^{\infty}$ implies the inverse of inequality \eq{2000}. This together with the monotonicity of $\Big(u_n\Big)_{n=0}^{\infty}$ yields the following inequality
\Eq{989}{
u_n-u_{n-1}\geq u_{n+1}-u_n\geq 0 \qquad (n\in\N).
}
This shows that $\Big(u_n-u_{n-1}\Big)_{n=1}^{\infty}$ is a decreasing sequence with non-negative terms. Using \prp{78}, we can conclude that $\Big(u_n-u_{n-1}\Big)_{n=1}^{\infty}$ is a subadditive sequence. This establishes our first assertion.\\

The decreasingness of the sequence $\Big(u_n-u_{n-1}\Big)_{n=1}^{\infty}$ can be denoted by the left-most inequality of \eq{989}. Arranging the terms we arrive at
\Eq{*}{
u_{n+1}+u_{n-1}\leq 2u_n \qquad (n\in\N).
}
It proves concavity of the sequence $\Big(u_n\Big)_{n=0}^{\infty}$ and establishes the result.
\end{proof}

These findings open several avenues for future research, including the characterization of various generalized forms of convex sequences. Similarly one can explore the stability results for newly derived subadditive sequences.
\bibliographystyle{plain}

\end{document}